\documentclass[11pt]{article}

\usepackage{amssymb, amsmath, amsthm, bbm, verbatim}

\theoremstyle{definition}
\newtheorem{defn}{Definition}

\newtheorem{lem}[defn]{Lemma}

\theoremstyle{plain}
\newtheorem{thm}[defn]{Theorem}

\newcommand{\bZ}{{\mathbb Z}}
\newcommand{\bQ}{{\mathbb Q}}

\newcommand{\bC}{{\mathbb C}}

\renewcommand{\r}{{\mathbb R}}
\newcommand{\z} {{\mathbb Z}}
\newcommand{\n} {{\mathbb N}}
\newcommand{\ddd}{,\dots,}
\renewcommand{\lll}{\left(}
\newcommand{\rrr}{\right)}

\newcommand{\supp}{\mbox{supp}\,}

\renewcommand{\a}{\alpha}

\newcommand{\be}{\begin{equation}}
\newcommand{\ee}{\end{equation}}
\newcommand{\ba}{\begin{eqnarray}}
\newcommand{\ea}{\end{eqnarray}}
\newcommand{\ban}{\begin{eqnarray*}}
\newcommand{\ean}{\end{eqnarray*}}

\begin{document}

\centerline{\Large\bf Haar bases for $L^2(\mathbb{Q}_2^2)$  generated by}
\centerline{\Large\bf one wavelet function
\footnote{The research was supported  by
DFG Project 436 RUS 113/951. The second author
also was supported  by Grant
 09-01-00162 of RFBR}
}

\vspace{.5cm}
\centerline{\bf S. Albeverio and M. Skopina}
\centerline{Inst. Appl. Math., HCM,
University of Bonn, e-mail: albeverio@uni-bonn.de}
\centerline{and}
\centerline{St. Petersburg State University, e-mail: skopina@MS1167.spb.edu}

\centerline{\bf Abstract}
The concept of
$p$-adic quincunx Haar MRA was introduced and studied in~\cite{KS10}.
In contrast to the real setting, infinitely many different wavelet
bases are generated by a $p$-adic MRA.
We give an explicit description for all wavelet functions
corresponding to the quincunx Haar MRA. Each one generates
an orthogonal basis, one of them was presented
in~\cite{KS10}. A connection between quincunx Haar bases and
two-dimensional separable Haar MRA is also found.



%%%%%%%%%%%%%%%%%%%%%%%%%
\section{Introduction}
%%%%%%%%%%%%%%%%%%%%%%%%%%%
In 2002, S.~V.~Kozyrev~\cite{Koz02} found a
$p$-adic wavelet basis for ${ L}^2(\Bbb Q_p)$ which is an analog of
the real Haar basis.    V.~M.~Shelkovich
 and one of the authors developed an  MRA approach
 to $p$-adic wavelets \cite{ShSk06}.  The scheme was
realized to construct the $p$-adic Haar MRA. In contrast to the real
setting, it appears  that there exist infinitly many
different orthonormal wavelet bases with the minimal
number of generating wavelet functions in the same Haar MRA. All these bases
were described explicitly in~\cite{ShSk06} for $p=2$
and in~\cite{Kh-Sh-Sk2} for arbitrary p,
one of these bases (for each $p$)  coincides with Kozyrev's wavelet basis.
 A wide class of orthogonal scaling functions generating
an MRA was described in~\cite{Kh-Sh-Sk1} by
A.~Khrennikov, V.~Shelkovuch and one of the authors.
However it has been shown in~\cite{AES2} by  S.~Evdokimov  and
 both  authors  that all these scaling functions
lead to the same Haar MRA and that there exist no other orthogonal
test scaling functions generating an MRA, except for those
described in~\cite{Kh-Sh-Sk1}. Thus all univariate orthogonal
MRA-based wavelet bases with a minimal possible number
of wavelet functions are described. It is not known if other
$p$-adic orthogonal wavelet bases exist.
Although, for the construction
  $p$-adic wavelets, there is another approach  based on $p$-adic wavelet set theory
introduced  by J.~J.~Benedetto and R.~L.~Benedetto~\cite{rlB04}, \cite{BB04},
 we did not see  up to now orthogonal $p$-adic
wavelet bases which are not generated by the Haar MRA.

The simplest way to construct a multivariate wavelet basis is the following.
Given an univariate MRA with a scaling function $\varphi$,
using the method suggested
by Y.~Meyer~\cite{M} (see, e.g.~\cite[Sec 2.1]{NPS}),
 one can easily construct a $d$-dimentional
 separable MRA,
where the scaling function is $\varphi\otimes\dots\otimes \varphi$.
If $\psi$ is a wavelet function in the univariate MRA, then the  functions
$f_1\otimes\dots\otimes f_d$, where each $f_k$ equals either $\varphi$ or
$\psi$ and not all $f_k$, $k=1,\dots, d$, are equal to $\varphi$,
form a set of wavelet functions generating a multivariate wavelet basis
(i.e. the basis consists of the dilations and shifts of the wavelet functions).
This approach was realized for $p$-adics  in~\cite{ShSk06}.
Next we are going to discuss only the case $p=d=2$. The separable
wavelet basis is generated by three wavalet functions in this case.
E. King and one of the authors studied two-dimensional
Haar MRA with quincunx dilation  in~\cite{KS10}. The  basis is generated by only
one wavelet function in this case. In the present paper we describe
all  wavelet functions corresponding to the quincunx MRA
and generating an orthogonal basis. We also
  discuss how these bases related to the separable MRA.

%%%%%%%%%%%%%%%%%%%%%%%%
\section{Preliminaries and notations}
%%%%%%%%%%%%%%%%%%%%%%%%%
%%%%%%%%%%%%%%%%%%%%%%%%
\subsection{$p$-adic numbers}

The theory of $p$-adic numbers is presented in detail in many books
(see, e.g.,\cite{Vl-V-Z}, \cite{Katok}).
We restrict ourselves to a brief description.

Let $\Bbb N$, $\Bbb Z$, $\Bbb R$, $\Bbb C$ be the sets of positive integers, integers,
real numbers, complex numbers, respectively.
The field $\Bbb Q_p$ of $p$-adic numbers is defined as the completion
of the field of rational numbers $\Bbb Q$ with respect to the
non-Archimedean $p$-adic norm $|\cdot|_p$. This $p$-adic norm
is defined as follows: $|0|_p=0$; if $x\ne 0$, $x=p^{\gamma}\frac{m}{n}$,
where $\gamma=\gamma(x)\in \Bbb Z$
and the integers $m$, $n$ are not divisible by $p$, then
$|x|_p=p^{-\gamma}$.
The norm $|\cdot|_p$ satisfies the strong triangle inequality
$|x+y|_p\le \max(|x|_p,|y|_p)$.
The canonical form of any $p$-adic number $x\ne 0$ is
\begin{equation}
\label{2}
x=p^{\gamma}\sum_{j=0}^\infty x_jp^j
%(x_0 + x_1p + x_2p^2 + \cdots),
\end{equation}
where $\gamma=\gamma(x)\in \Bbb Z$, \ $x_j\in D_p:=\{0,1,\dots,p-1\}$, $x_0\ne 0$.
%$j=0,1,\dots$.
The fractional part $\{x\}_p$ of the number $x$ equals by
definition $p^{\gamma}\sum_{j=0}^{-\gamma-1} x_jp^j$.
%:=p^{\gamma}(x_0 + x_1p  + \cdots+ x_{-\gamma-1}p^{-\gamma-1})$
Thus, $\{x\}_p=0$ if and only if $\gamma\ge0$. We also set $\{0\}_p=0$.

%%%%%%%%%%%%%%%%%%%%%%%%

For every prime $p$, the additive group of $\mathbb{Q}_{p}$
is a locally compact abelian group which contains the compact
open subgroup $\mathbb{Z}_{p}$.  The group $\mathbb{Z}_{p}$
is the set $\lbrace\alpha\in\mathbb{Q}_{p}\Big |\lvert\alpha\rvert_{p}\leq1\rbrace$,
the unit ball in $\mathbb{Q}_{p}$.  Equivalently, $\mathbb{Z}_p$
is the subgroup generated by $1$. It is well-known (\cite{RS00})
that since $\mathbb{Q}_p$ is a locally compact abelian group with a compact
open subgroup, it has a Haar measure normalized so that the measure of $\mathbb{Z}_p$
is $1$.  For simplicity, we shall denote this measure by $dx$.

 We also define the set
$
I_p = \{ x \in \mathbb{Q}_p : \{x\} = x \}.
$
$I_p$ is a set of coset representatives for $\mathbb{Q}_p / \mathbb{Z}_p$.
Since $\mathbb{Z}_p$ is open, $I_p$ is discrete.

For any $d \ge 1$, $\mathbb{Q}_p^d$ is a vector space over
$\mathbb{Q}_p$, i.e., $\Bbb Q_p^d$
consists of the vectors $x=(x_1,\dots,x_d)$,
where $x_j \in \Bbb Q_p$, $j=1,\dots,d$.
The $p$-adic norm on $\Bbb Q_p^d$ is
\begin{equation}
\nonumber
|x|_p:=\max_{1 \le j \le d}|x_j|_p.
\end{equation}

The Haar measure $dx$ on the field $\Bbb Q_p$ is extended to a
Haar measure $d^d x=dx_1\cdots dx_d$ on $\Bbb Q_p^d$ in the standard way.

Denote by $B_{N}=\{x\in \bQ^d_p: |x|_p \le p^{N}\}$
the ball of radius $p^{N}$ with the center at the point $0$,
$N \in \z$.

Finally, we define the set $I_p^d\subset \Bbb Q_p^d$ by
$
I_p^d:=I_p\times\dots\times I_p.
$

%%%%%%%%%%%%%%%%%%%%%%%%%%
\subsection{Functions of $p$-adic variables}

There are two main function theories over the $p$-adics.
One deals with functions $\mathbb{Q}_p^d \rightarrow \mathbb{C}$
and the other with $\mathbb{Q}_p^d \rightarrow \mathbb{Q}_p$.
We shall only deal with the former theory.

Since a Haar measure exists in $\Bbb Q_p^d$, the spaces ${ L}^q(\Bbb Q_p^d)$
can be naturally introduced.

The functions $\chi_p(x)=e^{2\pi i\{ x\}_p}$ are  additive characters for
the field $\bQ_p$.
 For any  set $S \subseteq \Bbb Q_p^d$,  $\mathbbm{1}_S$ denotes
  the characteristic function of $S$, i.e.
\begin{equation*}
\mathbbm{1}_S(x) = \left\{ \begin{array}{ccc} 1 & ; & x\in S \\ 0 &  ;
& \mathrm{else} \end{array} \right..
\end{equation*}
Note that $\mathbbm{1}_{\bZ_p}\otimes\mathbbm{1}_{\bZ_p}=\mathbbm{1}_{\bZ_p^2}$.

%%%%%%%%%%%%%%%%%%%%%%%%%

%%%%%%%%%%%%%%%%%%%%%%%
\subsection{$p$-Adic  MRA}\label{s6.1}

\begin{defn}[\cite{ShSk06}]
\label{de1}
A collection of closed spaces
$V_j\subset L^2(\bQ_p)$, $j\in\bZ$, is called a
{\it multiresolution analysis {\rm(}MRA{\rm)} in $L^2(\bQ_p)$} if the
following axioms hold

(a) $V_j\subset V_{j+1}$ for all $j\in\bZ$;

(b) $\bigcup\limits_{j\in\bZ}V_j$ is dense in $L^2(\bQ_p)$;

(c) $\bigcap\limits_{j\in\bZ}V_j=\{0\}$;

(d) $f(\cdot)\in V_j \Longleftrightarrow f(p^{-1}\cdot)\in V_{j+1}$
for all $j\in\bZ$;

(e) there exists a function $\varphi \in V_0$, called the {\em scaling function},
such that the system $\{\varphi(\cdot-a), a\in I_p\}$ is an orthonormal
basis for $V_0$.
\end{defn}

According to the standard scheme (see, e.g.,~\cite[\S~1.3]{NPS})
for the construction of MRA-based wavelets, for each $j$, we define
the space $W_j$ ({\em wavelet space}) as the orthogonal complement
of $V_j$ in $V_{j+1}$, i.e.,
\begin{equation}
\label{61}
V_{j+1}=V_j\oplus W_j, \qquad j\in \bZ,
\end{equation}
where $W_j\perp V_j$, $j\in \bZ$. It is not difficult to see that
\begin{equation}
\label{61.0}
f\in W_j \Longleftrightarrow f(p^{-1}\cdot)\in W_{j+1},
\quad\text{for all}\quad j\in \bZ
\end{equation}
and $W_j\perp W_k$, $j\ne k$.
Taking into account axioms (b) and (c), we obtain
\begin{equation}
\label{61.1}
{\bigoplus\limits_{j\in\bZ}W_j}=L^2(\bQ_p)
\quad \text{(orthogonal direct sum)}.
\end{equation}

If now we find   functions $\psi_\nu \in W_0$,
$\nu=1,\dots, r$,
such that the system $\{\psi_\nu(x-a), a\in~I_p, \nu=1,\dots, r\}$ is an orthonormal
basis for $W_0$, then, due to~(\ref{61.0}) and (\ref{61.1}),
the system $\{p^{j/2}\psi_\nu(p^{-j}\cdot-a), a\in I_p, j\in\bZ, \nu=1,\dots, r\}$,
is an orthonormal basis for $L^2(\bQ_p)$.
Such  functions $\psi_\nu$ are called  {\em wavelet functions} and
the basis is a {\em wavelet basis}.

\section{Separable and quincunx  $p$-adic wavelets bases}

\subsection{Separable  MRA }
Let $\{V_j^{(\nu)}\}_{j\in\bZ}$ , $\nu=1,\dots,n$, be one-dimensional
MRAs.
We introduce subspaces $V_j$, $j\in\bZ$, of $L^2(\bQ_p^d)$ by
\begin{equation}
\label{d3}
V_j=\bigotimes_{\nu=1}^dV^{(\nu)}_j=
\overline{{\rm span}\{F=f_1\otimes\dots\otimes f_n, \ f_\nu\in V_j^{(\nu)}\}}.
\end{equation}
Let $\varphi^{(\nu)}$ be a scaling function of $\nu$-th MRA $\{V_j^{(\nu)}\}_j$.
Set
$$
\Phi=\varphi^{(1)}\otimes\dots\otimes\varphi^{(n)}.
$$

It is not difficult to check that the following theorem holds (see~\cite{ShSk06}).
\begin{thm}
\label{th1}
Let $\{V_j^{(\nu)}\}_{j\in\bZ}$, $\nu=1,\dots,n$, be MRAs in $L^2(\bQ_p)$.
Then the subspaces $V_j$ of  $L^2(\bQ_p^d)$ defined by~{\rm(\ref{d3})}
satisfy the following properties:

{\rm(a)} $V_j\subset V_{j+1}$ for all $j\in\bZ$;

{\rm(b)} $\cup_{j\in\bZ}V_j$ is dense in $L^2(\bQ_p^d)$;

{\rm(c)} $\cap_{j\in\bZ}V_j=\{0\}$;

{\rm(d)} $f(\cdot)\in V_j \Longleftrightarrow f(p^{-1}\cdot)\in V_{j+1}$
for all $j\in\bZ$;

{\rm(e)} there exists a function $\Phi \in V_0$ such that the system
$\{\Phi(x-a), a\in I_p^d\}$, forms an orthonormal basis for $V_0$.
\end{thm}

The collection of spaces
$V_j$, $j\in\bZ$, discussed in Theorem~\ref{th1}
is called a {\it separable multiresolution analysis} in $L^2(\bQ_p^n)$,
the function $\Phi$ from axiom (e) is said to be its scaling function.

Following  the standard scheme (see, e.g.,~\cite[\S\  2.1]{NPS}),
we define the wavelet space $W_j$ as the orthogonal complement of $V_j$ in
$V_{j+1}$, i.e.,
$$
W_j=V_{j+1}\ominus V_j, \quad j\in\bZ.
$$
Since
$$
V_{j+1}=\bigotimes\limits_{\nu=1}^dV_{j+1}^{(\nu)}=
\bigotimes\limits_{\nu=1}^d\big(V^{(\nu)}_j\oplus W^{(\nu)}_j\big)
\qquad\qquad\qquad\qquad
$$
$$
\qquad
=V_j\oplus
\bigoplus\limits_{e\subset\{1,\dots,n\}, \, e\ne\emptyset}
\big(\bigotimes\limits_{\nu\in e} W^{(\nu)}_j \big)
\big(\bigotimes\limits_{\mu\not\in e} V^{(\mu)}_j \big).
$$
So, the space $W_j$ is a direct sum of $2^d-1$ subspaces
$W_{j,e}$,  $e\subset\{1,\dots,n\}$, $e\ne\emptyset$.

Let $\psi^{(\nu)}$ be a wavelet function,
i.e. a function whose shifts (with respect to
$a\in I_p$) form an orthonormal basis for $W^{(\nu)}_0$. It is clear
that the shifts (with respect to
$a\in I_p^d$) of the function
\begin{equation}
\label{wd-62}
\Psi_e=\big(\bigotimes\limits_{\nu\in e} \psi^{(\nu)}\big)
\big(\bigotimes\limits_{\mu\not\in e}\varphi^{(\mu)}\big),
\quad
e\subset\{1,\dots,n\}, \quad e\ne\emptyset,
\end{equation}
form an orthonormal basis for $W_{0,e}$.
So, we have
$$
L^2(\bQ_p^d)=\bigoplus\limits_{j\in\bZ}W_j
=\bigoplus\limits_{j\in\bZ}\Big(
\bigoplus\limits_{e\subset\{1,\dots,n\}, \, e\ne\emptyset}
W_{j,e}\Big),
$$
and  the  functions $p^{-dj/2}\Psi_e(p^j\cdot+a)$,
$e\subset\{1,\dots,n\}$, $e\ne\emptyset$, $j\in\bZ$, $a\in I_p^d$,
 form an orthonormal basis for $L^2(\bQ_p^d)$. A wavelet basis constructed in
 this way is called separable.

In particular, if   $p=d=2$, $\varphi^{(1)}=\varphi^{(2)}=\varphi$,
(which implies $\{V_j^{(1)}\}_{j\in\bZ}=\{V_j^{(2)}\}_{j\in\bZ}$)
and $\psi(x)=\chi_2(x/2)\mathbbm{1}_{\bZ_2}(x)$
is a wavelet function corresponding to the univariate  MRA,
then there are three wavelet functions
$
\Psi^1=\varphi\otimes\psi,\ \ \Psi^2=\psi\otimes\varphi,\ \ \Psi^3=\psi\otimes\psi.
$
corresponding to the separable MRA.

%%%%%%%%%%%%%%%%%%%%%%%%
\subsection{Quincunx  MRA}

Set
\begin{equation*}
A = \left( \begin{array}{cc} 1 & -1 \\ 1 & 1 \end{array} \right)^{-1} =
\left( \begin{array}{cc} \frac{1}{2} & \frac{1}{2} \\ -
\frac{1}{2} & \frac{1}{2} \end{array} \right).
\end{equation*}
$A$ is the inverse of the well-known quincunx matrix. In~\cite{GM92}
(see also~\cite[\S~2.8]{NPS}),
a Haar multiresolution analysis was presented for $L^2(\mathbb{R}^2)$
which used dilations by the quincunx.  Since $|\det A^{-1}| = 2$,
 there was only one wavelet generator.  However, the support
 of the scaling function was a fractal,
 the twin dragon fractal.  In contrast,
 it was proved in~\cite{KS10} that the Haar MRA
 for $L^2(\mathbb{Q}_2^2)$ associated to this matrix
 corresponds to a scaling function
 with a simple support, namely, $\mathbb{Z}_2^2$. To state this
 in more detail set
 \begin{equation}
 \label{99}
V_j^Q = \overline{\operatorname{span}
\{ \mathbbm{1}_{\bZ_2^2}(A^j \cdot - a): a \in I_2^2 \}}.
\end{equation}

\begin{thm}[\cite{KS10}]
\label{defn:pad_mra}
The subspaces $V^Q_j$ of  $L^2(\bQ_2^2)$, $j\in\z^d$, defined by~{(\ref{99})}
satisfy the following properties:
\begin{itemize}
\item[(a)] $V^Q_j \subset V^Q_{j+1}$ for all $j \in \mathbb{Z}$,
\item[(b)] $\bigcup_{j\in\mathbb{Z}} V^Q_j$ is dense in $L^2(\mathbb{Q}^2_2)$,
\item[c()] $\bigcap_{j \in \mathbb{Z}} V^Q_j = \{0\}$,
\item[(d)] $f \in V^Q_j \Leftrightarrow f(A\cdot)
\in V^Q_{j+1}$ for all $j\in \mathbb{Z}$,
\item[(e)] there exists a function $\phi \in V^Q_0$,
%called the \emph{scaling function},
such that the system
$\{\phi (\cdot - a) : a \in I_2^2 \}$ is an orthonormal basis for $V^Q_0$.
\end{itemize}
\end{thm}

The collection of spaces
$V^Q_j$, $j\in\bZ$, discussed in Theorem~\ref{defn:pad_mra}
is called the {\it quincunx Haar multiresolution analysis} in $L^2(\bQ_p^n)$,
the function $\phi$ from axiom (e) is its scaling function.
It follows from~(\ref{99}) that $\phi=\mathbbm{1}_{\bZ_2^2}$.

It was proved in~\cite{KS10} that the function
\be
\psi = \phi(A\cdot) - \phi\lll A\cdot - \binom{1/2}{1/2}\rrr
\label{97}
\ee
is a wavelet function corresponding to the quincunx MRA,
i.e. the system
$\{\psi(\cdot-a), a\in I_p^2\}$ is an orthonormal basis for
the wavelet space $W^Q_0:=V^Q_1\ominus V^Q_0$. It follows that
the system
\be
\{\psi(A^j\cdot-a), \ \ j\in\z, a\in I_2^2\}
\label{98}
\ee
is an orthonormal basis for $L^2(\bQ_2^2)$. %This basis is called
%{\em quincunx wavelet bases}

Next we shall show that this basis corresponds also to the
separable MRA in some sense. Let $\{V_j\}_{j\in\bZ}$ be a
separable MRA with  $p=d=2$, $\varphi^{(1)}=\varphi^{(2)}=\mathbbm{1}_{\bZ_2}$,
i.e. both $\{V_j^{(1)}\}_{j\in\bZ}$ and  $\{V_j^{(2)}\}_{j\in\bZ}$
are  Haar MRAs. In this case, $\Phi=\mathbbm{1}_{\bZ_2}\otimes\mathbbm{1}_{\bZ_2}=
\mathbbm{1}_{\bZ_2^2}$.
Thus, $\Phi=\phi$, which yields $V_0=V^Q_0$. Since
\begin{equation*}
A^2 = \left( \begin{array}{cc} 0 & 1/2 \\ 1/2 & 0 \end{array} \right),
\end{equation*}
we have
$$
\phi(A^2x-a)=\phi\lll\binom{\frac{x_2}2-a_1}{\frac{x_1}2-a_2}\rrr=
\phi\lll\binom{\frac {x_1}2-a_2}{\frac{x_2}2-a_1}\rrr.
$$
It follows that
$$
\{\phi(A^2x-a), \ \ a\in I_2^2\}=\{\phi(\frac x2-a), \ \ a\in I_2^2\}.
$$
So,
$
V^Q_0\oplus W^Q_0\oplus W^Q_1=V^Q_2=V_1,
$
which yields that
$
W_0=W^Q_0\oplus W^Q_1,
$
the functions
$
 \psi(\cdot-a), \psi(A\cdot-a), \ \ a\in I_2^2,
$
form an orthonormal basis for $W_0$ and the dyadic dilations
of these functions
form an orthonormal basis for $L^2(\bQ_2^2)$. These relations
may be treated as follows: the quincunx Haar MRA is the ''square root''
of the separable Haar MRA.
%%%%%%%%%%%%%%%%%%%%%%%%%%

\section{Description of wavelet functions}

As was said above, (\ref{98})
is an orthonormal basis for $L^2(\mathbb{Q}_2^2)$ generated by a
single wavelet function $\psi$. In contrast to the real setting,
 $p$-adic wavelet function corresponding to an MRA is not unique.
Now we are going to find other wavelet functions corresponding to the
quincunx Haar MRA.

%%%%%%%%%%%%%%%%%%%%

First we shall prove some  auxiliary facts.
\begin{lem}
Let $\psi$, $\tilde\psi$ be compactly supported functions with
their supports contained in $B_s$, $s\in\n$, and such that
\be
\tilde\psi(x)=\sum_{a\in I_2^2}h_a\psi(x-a).
\label{96}
\ee
Then $h_a=0$ whenever $|a|_2>2^s$.
\label{l1}
\end{lem}
{\bf Proof.} We rewrite~(\ref{96}) in the form
$$
\tilde\psi(x)=
\sum_{a\in I^2_2\atop|a|_2\le2^s}h_a\psi(x-a)+
\sum_{a\in I^2_2\atop|a|_2>2^s}h_a\psi(x-a).
$$
The second sum on the right-hand side vanishes whenever
$|x|_2\le2^s$ because $\supp\psi\subset B_s$.
The first sum on the right-hand side and the left-hand side
 vanish whenever $|x|_2>2^s$ because $\supp\psi\subset B_s$
and $\supp\tilde\psi\subset B_s$ respectively.
So, the second sum is  equal to zero for all $x\in \bQ_2^2$. $\Diamond$

\begin{lem}
Let $n\in\z_2^2$, and let $\psi$ be the function defined in~(\ref{97}).
If $An\in\z_2^2$, then $\psi(\cdot+n)=\psi$,
if $An\not\in\z_2^2$, then $\psi(\cdot+n)=-\psi$.
\label{l2}
\end{lem}
{\bf Proof.} Observe that $An\in\z_2^2$ if and only if
both coordinates of $n$ are either even or odd. Therefore,
if $An\in\z_2^2$, then both coordinates of $An$ are integers,
otherwise both coordinates of $An$ are semi-integers.
It remains to note that the function $\phi$ is $1-periodic$
with respect to each variable. $\Diamond$

\begin{thm}
All compactly supported  wavelet functions corresponding to the
quincunx Haar MRA  are given by
\be
\tilde\psi(x)=\sum\limits_{k=0}^{2^s-1}\sum\limits_{l=0}^{2^s-1}
\alpha_{kl}\psi\lll x-\binom{k/2^s}{l/2^s}\rrr
\label{101}
\ee
where $s\in\n$, $\psi$ is the function defined by~(\ref{97}),
\be
\alpha_{kl}=\left\{
\begin{array}{ll}
2^{-2s}(-1)^ke^{-\frac{\pi ik}{2^s}}
\sum\limits_{p=0}^{2^s-1}\sum\limits_{q=0}^{2^s-1}
e^{-2\pi i\frac{qk}{2^s}}\gamma_{qp},
&\mbox{if\ \ \ } l=0,
\\
2^{-2s}(-1)^{k-l+1}e^{-\frac{\pi i(k-l)}{2^s}}
\sum\limits_{p=0}^{2^s-1}\sum\limits_{q=0}^{2^s-1}
e^{-2\pi i\frac{qk-lp}{2^s}}\gamma_{qp},
&\mbox{if\ \ \ } l\ne0,
\end{array}
\right.
\label{95}
\ee
$$
\gamma_{qp}\in\bC,\ \ |\gamma_{qp}|=1,\ \ p,q=0\ddd2^s-1.
$$
\label{maintheorem}
\end{thm}

{\bf Proof.} Let $\tilde\psi$ be a compactly supported wavelet
function corresponding to the quincunx Haar MRA, i.e., the
system  $\{\tilde\psi(\cdot-a), \ \  a\in I_2^2\}$ is an orthonormal
basis for $W_0^Q$, and $\supp\tilde\psi\subset B_s$. Since
$\{\psi(\cdot-a), \ \  a\in I_2^2\}$ is a basis for $W^Q_0$
and $\tilde\psi\in W^Q_0$, we have
$$
\tilde\psi(x)=\sum_{a\in I_2^2}h_a\psi(x-a).
$$
Due to Lemma~\ref{l1}, taking into account that
$\supp\psi\subset B_1\subset B_s$, $\supp\tilde\psi\subset B_s$,
we obtain~(\ref{101})
with some coefficients $\alpha_{kl}$.

Since $\{\psi(\cdot-a), \ \  a\in I_2^2\}$ is an orthogonal system,
evidently, $\tilde\psi$ is orthogonal to $\tilde\psi(\cdot-a)$ whenever
$a\in I_2^2$ and $a\ne\binom{k/2^s}{l/2^s}, k,l=0,1\ddd 2^s-1$.
Thus, the system $\{\tilde\psi(\cdot-a), \ \  a\in I_2^2\}$ is an orthonormal
system if and only if the system consisting of the functions
$\tilde\psi\lll x-\binom{k/2^s}{l/2^s}\rrr$, $k,l=0,1\ddd 2^s-1$, is orthonormal.

Next we need the following notations. Set
$$
\Lambda=\left( \begin{array}{ccccc}
\mathbbm{A} & \Bbb O&\hdots&\Bbb O&\Bbb O \\
\Bbb O & \mathbbm{A} &\hdots&\Bbb O&\Bbb O\\
\vdots & \vdots &\ddots&\vdots&\vdots\\
\Bbb O & \Bbb O &\hdots&\mathbbm{A}&\Bbb O\\
\Bbb O & \Bbb O &\hdots&\Bbb O&\mathbbm{A}
\end{array} \right),
\ \ \ \
\Omega=\left( \begin{array}{ccccc}
\Bbb O & \Bbb O&\hdots&\Bbb O& -\mathbbm{I}\\
\mathbbm{I} & \Bbb O &\hdots&\Bbb O&\Bbb O\\
\Bbb O & \mathbbm{I} &\hdots&\Bbb O&\Bbb O\\
\vdots & \vdots &\ddots&\vdots&\vdots\\
\Bbb O & \Bbb O &\hdots&\mathbbm{I}&\Bbb O
\end{array} \right),
$$
where $\Bbb O,\mathbbm{I}, \mathbbm{A}$ are $2^s\times 2^s$ matrices,
$\Bbb O$ is the zero matrix, $\mathbbm{I}$ is the unity matrix,
$$
\mathbbm{A}=\left( \begin{array}{ccccc}
0 & 0&\hdots&0& -1\\
1 & 0 &\hdots&0&0\\
0 & 1 &\hdots&0&0\\
\vdots & \vdots &\ddots&\vdots&\vdots\\
0 & 0 &\hdots&1&0
\end{array} \right).
$$
Let $\a$ be a $2^{2s}$-dimensional vector (column) whose
$N$-th coordinate, $N=2^sl+k$, $k,l=0,1\dots 2^s-1$,
is $\a_{kl}$, and let
$\Psi$, $\tilde\Psi$ be $2^{2s}$-dimensional vector-functions  (columns) whose
$N$-th coordinates, $N=2^sl+k$, $k,l=0,1\ddd 2^s-1$, are
$\psi\lll x-\binom{k/2^s}{l/2^s}\rrr$ and
$\tilde\psi\lll x-\binom{k/2^s}{l/2^s}\rrr$ respectively.
Using Lemma~\ref{l2}, we obtain $\tilde\Psi=D\Psi$, where
$D$ is a $2^{2s}\times2^{2s}$ matrix whose $N$-th
row, $N=2^sl+k$, $k,l=0,1\ddd 2^s-1$, is $\Omega^l\Lambda^ka$.

Due to orthonormality of the system $\{\psi(\cdot-a), \ \  a\in I_2^2\}$,
the coordinates of $\tilde\Psi$ form an orthonormal system
if and only if the matrix $D$ is unitary. To describe all unitary
matrices $D$ we need to find all  $\a\in {\Bbb C}^{2^{2s}}$
such that the vectors
$\a, \Lambda \a\ddd\Omega^l\Lambda^k \a\ddd\Omega^{2^s-1}\Lambda^{2^s-1}\a$
are orthonormal. We have already one such  vector $\a_0=(1,0\ddd0)^T$
because the matrix $D_0$ whose rows are
$\a_0, \Lambda \a_0\ddd\Omega^l\Lambda^k \a_0\ddd\Omega^{2^s-1}\Lambda^{2^s-1}\a_0$
is, evidently, unitary.

Let us prove that the rows of $D$ are orthonormal if and only if
$\a=B\a_0$, where $B$ is a unitary matrix such that
\be
\Lambda B=B\Lambda,
\label{104}
\ee
\be
\Omega B=B\Omega.
\label{105}
\ee
Indeed, let $\a=B\a_0$, with $B$ being a unitary matrix satisfying~(\ref{104}),
(\ref{105}). In this case
$
\Omega^l\Lambda^k\a=\Omega^l\Lambda^kB\a_0=B(\Omega^l\Lambda^k\a_0),
$
and the unitarity of $D_0$ implies the unitarity of $D$.
Conversly, if the rows of $D$ are orthonormal, taking into account that
the rows of $D_0$ are also orthonormal, we conclude that
there exists a unitary matrix $B$ such that
\be
D^T=BD_0^T
\label{106}
\ee
So we have, in particular,
$$
\a=B\a_0,\ \  \Lambda\a=B\Lambda\a_0,\ \ \Lambda^{2}\a=B\Lambda^{2}\a_0
 \ddd\Lambda^{2^s-1}\a=B\Lambda^{2^s-1}\a_0.
$$
Substituting the first equality into the second one, we obtain
$(\Lambda B-B\Lambda)\a_0=0$. Using this and substituting
the first equality into the third one,  we obtain
$(\Lambda B-B\Lambda)(\Lambda \a_0)=0$. After  similar manipulations
with other equalities we have
\be
(\Lambda B-B\Lambda)(\Lambda^p \a_0)=0,\ \ p=0,1\ddd 2^s-2.
\label{107}
\ee
Since $\Lambda^{2^s}\a=-\a$, $\Lambda^{2^s}\a_0=-\a_0$, substituting
$\a=B\a_0$ and using~(\ref{107}) we obtain
$(\Lambda B-B\Lambda)(\Lambda^{2^s-1} \a_0)=0$.
Next it follows from~(\ref{106}) that
$$
\Omega\a=B\a_0,\ \  \Omega\Lambda\a=B\Omega\Lambda\a_0,\ \
\Omega\Lambda^{2}\a=B\Omega\Lambda^{2}\a_0
 \ddd\Omega\Lambda^{2^s-1}\a=B\Omega\Lambda^{2^s-1}\a_0.
$$
Substituting $\a=B\a_0$ into the first equality, we have
$(\Omega B-B\Omega)\a_0=0$. Using this, taking into account
that $\Omega\Lambda=\Lambda\Omega$ and substituting $\a=B\a_0$ into the
second equality, we obtain that $(\Lambda B-B\Lambda)(\Omega\a_0)=0$.
After similar manipulations with the other equalities and with the
equalities $\Omega\Lambda^{2^s}\a=-\Omega\a$,
$\Omega\Lambda^{2^s}\a_0=-\Omega\a_0$, we have
$$
(\Lambda B-B\Lambda)(\Lambda^p \Omega^q\a_0)=0,\ \ p=0,1\ddd 2^s-1.
$$
Continuing this process, we get that
$(\Lambda B-B\Lambda)(\Lambda^p \a_0)=0$ for all $p,q=0,1\ddd2^s-1$.
Since  $\{\Lambda^p \Omega^q\a_0\}_{p,q=0}^{2^s-1}$ is a basis
for $\bC^{2^{2s}}$, we obtain~(\ref{104}). In a similar way we can
check that (\ref{105}) holds.

Thus we should describe all unitary matrices $B$
satisfying~(\ref{104}), (\ref{105}).
Let
$$
B=\left( \begin{array}{ccc}
B_{00} & \hdots&B_{0,2^s-1}\\
 \vdots &\ddots&\vdots\\
B_{2^s-1,0} &\hdots&B_{2^s-1,2^s-1}
\end{array} \right),
$$
where $B_{ij}$ is a $2^s\times2^s$ matrix.
It follows from~(\ref{105}) that
{\footnotesize
$$
\left( \begin{array}{cccc}
-B_{2^s-1,0}&-B_{2^s-1,1} & \hdots&-B_{2^s-1,2^s-1}\\
B_{00} &B_{01}& \hdots&B_{0,2^s-1}\\
\vdots& \vdots &\ddots&\vdots\\
B_{2^s-2,0} &B_{2^s-2,1}&\hdots&B_{2^s-2,2^s-1}
\end{array} \right)
=
\left( \begin{array}{cccc}
B_{01} & \hdots&B_{0,2^s-1}&-B_{00}\\
B_{11} & \hdots&B_{1,2^s-1}&-B_{10}\\
 \vdots &\ddots&\vdots&\vdots\\
B_{2^s-1,1} &\hdots&B_{2^s-1,2^s-1}&-B_{2^s-1,0}
\end{array} \right).
$$
}
This yields
\ban
&&B_{00}=B_{11}=\dots=B_{2^s-1,2^s-1},
\\
&&B_{0l}=B_{1,l+1}=\dots=B_{2^s-l-1,2^s-1}=
\\
&&\hspace{3cm}-B_{2^s-l,0}=\dots=-B_{2^s-1,l-1},
\ \ l=1\ddd 2^s-1.
\ean
i.e.
$$
B=
\left( \begin{array}{cccc}
\beta_{0}&\beta_{1} & \hdots&\beta_{2^s-1}\\
-\beta_{2^s-1} &\beta_{0}& \hdots&\beta_{2^s-2}\\
\vdots& \vdots &\ddots&\vdots\\
-\beta_1 &-\beta_2&\hdots&\beta_{0}
\end{array} \right),
$$
where $\beta_l=B_{0l}$. It is not difficult to see that
any such matrix $B$ satisfies~(\ref{105}).

Since
{\footnotesize
$$
\Lambda B=\left( \begin{array}{cccc}
{\Bbb A}\beta_{0}&{\Bbb A}\beta_{1} & \hdots&{\Bbb A}\beta_{2^s-1}\\
{\Bbb A}(-\beta_{2^s-1}) &{\Bbb A}\beta_{0}& \hdots&{\Bbb A}\beta_{2^s-2}\\
\vdots& \vdots &\ddots&\vdots\\
{\Bbb A}(-\beta_1) &{\Bbb A}(-\beta_2)&\hdots&{\Bbb A}\beta_{0}
\end{array} \right),
\ \ \
B\Lambda=\left( \begin{array}{cccc}
\beta_{0}{\Bbb A}&\beta_{1}{\Bbb A} & \hdots&\beta_{2^s-1}{\Bbb A}\\
(-\beta_{2^s-1} ){\Bbb A}&\beta_{0}{\Bbb A}& \hdots&\beta_{2^s-2}{\Bbb A}\\
\vdots& \vdots &\ddots&\vdots\\
(-\beta_1){\Bbb A} &(-\beta_2){\Bbb A}&\hdots&\beta_{0}{\Bbb A},
\end{array} \right)
$$
}
(\ref{104}) holds if and only if $\beta_\nu{\Bbb A}={\Bbb A}\beta_\nu$,
for all $\nu=0\ddd2^s-1$. It follows from the proof of Theorem 4.1
in~\cite{ShSk06} that $\beta_n$ commutes with ${\Bbb A}$ if and only if
$\beta_\nu=C\tilde\beta_\nu C^{-1}$, where $\tilde\beta_\nu$ is a diagonal
matrix, $C$ is a matrix whose entries  are given by
\be
c_{pq}=2^{-s/2}(-1)^pe^{-\pi ip\frac{2q+1}{2^s}},\ \ \  p,q=0\ddd2^s-1.
\label{108}
\ee
Since $C$ is a unitary matrix, $B$ is unitary if and only if
the matrix
$$
\tilde B=
\left( \begin{array}{cccc}
\tilde\beta_{0}&\tilde\beta_{1} & \hdots&\tilde\beta_{2^s-1}\\
-\tilde\beta_{2^s-1} &\tilde\beta_{0}& \hdots&\tilde\beta_{2^s-2}\\
\vdots& \vdots &\ddots&\vdots\\
-\tilde\beta_1 &-\tilde\beta_2&\hdots&\tilde\beta_{0}
\end{array} \right)
$$
is unitary.
Set $\theta_l=(\lambda_l^{(0)}\ddd \lambda_l^{(2^s-1)})^T$, where
$\lambda_l^{(\nu)}$ is the $l$-th diagonal element of $\tilde\beta_\nu$.
It is easy to see that $\tilde B$ is unitary if and only if
the vectors $\theta_l, {\Bbb A}\theta_l\ddd   {\Bbb A}^{2^s-1}\theta_l$
form an orthonormal system for all $l=0\ddd2^s-1$.
It follows from the proof of Theorem 4.1 in~\cite{ShSk06} that all
such $\theta_l$ are given by
$$
\theta_l=C
\left( \begin{array}{cccc}
\gamma_{l0}&0 & \hdots&0\\
0 & \gamma_{l1}&\hdots&0\\
\vdots& \vdots &\ddots&\vdots\\
0&0&\hdots&\gamma_{l, 2^s-1}
\end{array} \right)
C^{-1}e_1,
$$
where $\gamma_{lk}\in\bC$, $|\gamma_{lk}|=1$,
$e_1=(1,0\ddd 0)^T\in\r^{2^s}$. Substituting~(\ref{108}),
we have
\be
\lambda_l^{(\nu)}=2^{-s}(-1)^\nu e^{-\pi i\frac{\nu}{2^s}}
\sum\limits_{k=0}^{2^s-1}e^{-2\pi i\frac{k\nu}{2^s}}\gamma_{lk},
\ \ \ l,\nu=0\ddd2^s-1.
\label{109}
\ee
Again using~(\ref{108}), we obtain
\ba
\a=B\a_0=
\left( \begin{array}{cccc}
C&{\Bbb O} & \hdots&{\Bbb O}\\
{\Bbb O} & C&\hdots&{\Bbb O}\\
\vdots& \vdots &\ddots&\vdots\\
{\Bbb O}&{\Bbb O}&\hdots&C
\end{array} \right)
\tilde B
\left( \begin{array}{cccc}
C^{-1}&{\Bbb O} & \hdots&{\Bbb O}\\
{\Bbb O} & C^{-1}&\hdots&{\Bbb O}\\
\vdots& \vdots &\ddots&\vdots\\
{\Bbb O}&{\Bbb O}&\hdots&C^{-1}
\end{array} \right)
\a_0=
\nonumber
\\
\left( \begin{array}{cccc}
C&{\Bbb O} & \hdots&{\Bbb O}\\
{\Bbb O} & C&\hdots&{\Bbb O}\\
\vdots& \vdots &\ddots&\vdots\\
{\Bbb O}&{\Bbb O}&\hdots&C
\end{array} \right)
\left( \begin{array}{c}
\delta_0\\
-\delta_{2^s-1}\\
\vdots\\
-\delta_1
\end{array} \right),
\label{110}
\ea
where
\be
\delta_\nu=\tilde\beta_\nu
\left( \begin{array}{c}
\overline{c_{01}}\\
\vdots\\
\overline{c_{0,2^s-1}}
\end{array} \right)=
2^{-s/2}\tilde\beta_\nu
\left( \begin{array}{c}
1\\
\vdots\\
1
\end{array} \right)=
2^{-s/2}
\left( \begin{array}{c}
\lambda_0^{(\nu)}\\
\vdots\\
\lambda_{2^s-1}^{(\nu)}
\end{array} \right).
\label{111}
\ee
To prove~(\ref{95}) it remains to combine~(\ref{110})
with~(\ref{111}), (\ref{109}), and take into account that
$\a_{0k}$ is the $k$-th coordinate of $C\delta_0$ and
$\a_{lk}$ is the $k$-th coordinate of $C\delta_{2^s-l}$
for $l=1\ddd 2^s-1$.\ $\Diamond$

%%%%%%%%%%%%%%%%%%%%%%%%%

\end{document}